\documentclass[12pt,psamsfonts,leqno,oneside,letterpaper]{amsart}
\usepackage[dvips,text={6.5truein,9truein},left=1truein,top=1truein]{geometry}
\usepackage{amssymb,amsmath,amscd,enumerate}
\usepackage[pdftex]{graphicx}
\usepackage[colorlinks,linkcolor=blue,citecolor=blue,pdfstartview=FitH]{hyperref}
\input xy
\xyoption{all}
\SelectTips{cm}{12}

\theoremstyle{definition}
\newtheorem{para}{}[section]

\newtheorem{remark}[para]{Remark}
\newtheorem{remarks}[para]{Remarks}
\newtheorem{notation}[para]{Notation}
\newtheorem{convention}[para]{Convention}
\newtheorem{definition}[para]{Definition}
\newtheorem{definitions}[para]{Definitions}

\newcommand\Alternatives{\begin{enumerate}[(i)]}
\newcommand\EndAlternatives{\end{enumerate}}
\newcommand\Conditions{\begin{enumerate}[(1)]}
\newcommand\EndConditions{\end{enumerate}}

\theoremstyle{plain}
\newtheorem{theorem}[para]{Theorem}
\newtheorem{lemma}[para]{Lemma}
\newtheorem{proposition}[para]{Proposition}
\newtheorem{corollary}[para]{Corollary}
\newtheorem{conjecture}[para]{Conjecture}
\newtheorem{claim}[equation]{}

\numberwithin{equation}{para}
\numberwithin{figure}{section}

\newcommand\Number{\begin{para}}
\newcommand\EndNumber{\end{para}}
\newcommand\Definition{\begin{definition}}
\newcommand\EndDefinition{\end{definition}}
\newcommand\Definitions{\begin{definitions}}
\newcommand\EndDefinitions{\end{definitions}}
\newcommand\Theorem{\begin{theorem}}
\newcommand\EndTheorem{\end{theorem}}
\newcommand\Conjecture{\begin{conjecture}}
\newcommand\EndConjecture{\end{conjecture}}
\newcommand\Remark{\begin{remark}}
\newcommand\EndRemark{\end{remark}}
\newcommand\Remarks{\begin{remarks}}
\newcommand\EndRemarks{\end{remarks}}
\newcommand\Convention{\begin{convention}}
\newcommand\EndConvention{\end{convention}}
\newcommand\Notation{\begin{notation}}
\newcommand\EndNotation{\end{notation}}
\newcommand\Lemma{\begin{lemma}}
\newcommand\EndLemma{\end{lemma}}
\newcommand\Proposition{\begin{proposition}}
\newcommand\EndProposition{\end{proposition}}
\newcommand\Corollary{\begin{corollary}}
\newcommand\EndCorollary{\end{corollary}}
\newcommand\Claim{\begin{claim}}
\newcommand\EndClaim{\end{claim}}
\newcommand\Proof{\begin{proof}}
\newcommand\EndProof{\end{proof}}
\newcommand\Equation{\begin{equation}}
\newcommand\EndEquation{\end{equation}}
\newcommand\NoProof{{\hfill$\square$}}
\newcommand\Bullets{\begin{itemize}}
\newcommand\EndBullets{\end{itemize}}

\newcommand\homologystuff{3.5}
\newcommand\topeleven{8.13}
\newcommand\topsix{1.1}
\newcommand\threefreevolume{9.3}

\newcommand\ifnotwhynot{7.1}

\newcommand\geomsix{1.2}

\newcommand\tX{\widetilde X}
\newcommand{\cut}{\,\backslash\backslash\,}

\newcommand\RR{{\bf R}}

\newcommand\chibar{\bar\chi}

\newcommand\inter{\mathop{\rm int}}

\newcommand\tM{\widetilde M}

\newcommand\tP{\widetilde P}

\newcommand\calb{{\mathcal B}}

\newcommand\cald{{\mathcal D}}

\newcommand\calp{{\mathcal P}}

\newcommand\calw{{\mathcal W}}

\newcommand\kish{\mathop{\rm kish}}

\newcommand\sosmall{shallow}
\newcommand\Hg{{\rm Hg}}
\newcommand\pile{layering}

\newcommand\ZZ{{\mathbb Z}}
\newcommand\Z{{\mathbb Z}}

\newcommand\HH{{\mathbb H}}

\newcommand\cala{{\mathcal A}}

\newcommand\cals{{\mathcal S}}

\newcommand\vol{\mathop{\rm vol}}

\begin{document}

\title{Incompressible surfaces, hyperbolic volume, Heegaard genus 
and homology}

\author{Marc Culler}
\address{Department of Mathematics, Statistics, and Computer Science
(M/C 249)\\
University of Illinois at Chicago\\
851 S. Morgan St.\\
Chicago, IL 60607-7045}
\email{culler@math.uic.edu}
\thanks{Partially supported by NSF grants DMS-0204142 and DMS-0504975}

\author{Jason DeBlois}
\address{Department of Mathematics, Statistics, and Computer Science
(M/C 249)\\
University of Illinois at Chicago\\
851 S. Morgan St.\\
Chicago, IL 60607-7045}
\email{jdeblois@math.uic.edu}
\thanks{Partially supported by NSF grant DMS-0703749}

\author{Peter B.~Shalen}
\address{Department of Mathematics, Statistics, and Computer Science
(M/C 249)\\
University of Illinois at Chicago\\
851 S. Morgan St.\\
Chicago, IL 60607-7045}
\email{shalen@math.uic.edu}
\thanks{Partially supported by NSF grants DMS-0204142 and DMS-0504975}

\begin{abstract}
We show that if $M$ is a complete, finite--volume, hyperbolic
$3$-manifold having exactly one cusp, and if
$\mathrm{dim}_{\mathbb{Z}_2} H_1(M;\mathbb{Z}_{2}) \geq 6$, then $M$
has volume greater than $5.06$.  We also show that if $M$ is a closed,
orientable hyperbolic 3--manifold with $\mathrm{dim}_{\mathbb{Z}_2}
H_1(M; \mathbb{Z}_2) \geq 4$, and if the image of the cup product map
$H^1 (M; \mathbb{Z}_{2})\otimes H^1 (M; \mathbb{Z}_{2})\to H^2 (M;
\mathbb{Z}_2)$ has dimension at most $1$, then $M$ has volume greater
than $3.08$. The proofs of these geometric results involve new
topological results relating the Heegaard genus of a closed Haken
manifold $M$ to the Euler characteristic of the kishkes of the
complement of an incompressible surface in $M$.
\end{abstract}

\maketitle

\section{Introduction}

If $S$ is a properly embedded surface in a compact 3-manifold $M$, let
$M \cut S$ denote the manifold which is obtained by cutting along $S$;
it is homeomorphic to the complement in $M$ of an open regular
neighborhood of $S$.

The topological theme of this paper is that the bounded manifold
obtained by cutting a topologically complex closed simple Haken
3-manifold along a suitably chosen incompressible surface $S\subset M$
will also be topologically complex.  Here the
``complexity'' of $M$ is measured by its Heegaard genus, and the
``complexity'' of $M\cut S$ is measured by the absolute value of the
Euler characteristic of its ``kishkes'' (see Definitions \ref{simple
  def} below).  

Our topological theorems have geometric consequences illustrating a
longstanding theme in the study of hyperbolic 3-manifolds --- that the
volume of a hyperbolic 3-manifold reflects its topological complexity.
We obtain lower bounds for volumes of closed and one-cusped hyperbolic
manifolds with sufficient topological complexity, extending work of
Culler and Shalen along the same lines. Here ``topological
complexity'' is measured in terms of the mod $2$ first homology, or
the mod $2$ cohomology ring.

\Definitions\label{simple def} We shall say that a $3$-manifold $M$ is
{\it simple} if
\begin{enumerate}[(i)]
\item $M$ is compact, connected, orientable, irreducible and
  boundary-irreducible;
\item no subgroup of $\pi_1(M)$ is isomorphic to $\Z\times\Z$; and
\item $M$ is not a closed manifold with finite fundamental group.
\end{enumerate}

Let $X$ be a simple $3$-manifold with $\partial X\ne\emptyset$.
According to \cite{JaS} or
\cite{Jo}, the characteristic submanifold $\Sigma_X$ of $X$ is
well-defined up to isotopy, and each component of $\Sigma_X$ is either
an $I$-bundle meeting $\partial X$ in its associated $\partial
I$-bundle, or a solid torus meeting $\partial X$ in a collection of
disjoint annuli that are homotopically non-trivial in $X$.  We define
$\kish(X)$ (the ``kishkes'' of $X$,  or ``guts'' in the terminology of
\cite{ASTD}) to be the union of all
components of $\overline{X-\Sigma_X}$ that have negative Euler
characteristic. The components of the frontier of $\kish(X)$ are
essential annuli in $X$.

If $X$ is a compact $3$-manifold whose components are all bounded and
simple and if $X_1,\ldots,X_k$ denote the components of $X$, we define
$\kish(X)=\kish(X_1)\cup\cdots\cup\kish(X_k)\subset X$.
\EndDefinitions

\Definition Let $g$ be an integer $\ge2$, let $h$ be a positive real
number, and let $M$ be an orientable, irreducible $3$-manifold. We
shall say that $M$ is {\it $(g,h)$-small} if every connected
  closed incompressible surface in $M$ has genus at least $h$ and
every separating connected closed incompressible surface in
$M$ has genus at least $g$.
\EndDefinition

We shall denote the Heegaard genus of a $3$-manifold $Q$ by
$\Hg(Q)$.

\newtheorem*{trichotomyThm}{Theorem \ref{trichotomy}}
\newcommand\trichotomy{ 
  Suppose $M$ is a closed, simple 3-manifold
  containing a separating connected closed incompressible
  surface of some genus $g$, that $\Hg(M) \geq g+4$, and
  that $M$ is $(g,\frac{g}{2}+1)$-small. Then $M$ contains a
  separating connected closed incompressible surface $S$ of
  genus $g$ satisfying at least one of the following conditions:
  \begin{enumerate}
  \item  at least one component of $M\cut S$ is acylindrical;  or
  \item for each  component $B$ of $M \cut S$ we have
  $\kish(B)\ne\emptyset$. 
 \end{enumerate}
}
\begin{trichotomyThm}
\trichotomy
\end{trichotomyThm}

The key idea in the proof is an organizing principle for cylinders
properly embedded in the complement of a separating connected closed
incompressible surface.  This is discussed in Sections
\ref{sec:annulus} and \ref{sec:trichotomy}.  We apply Theorem
\ref{trichotomy} in conjuction with the theorem below concerning
nonseparating surfaces, which is proved in Section \ref{sec:nonsep}.
For a manifold $M$ with (possibly empty) boundary, let $\chi(M)$
denote the Euler characteristic of $M$, and let $\bar{\chi}(M) =
-\chi(M)$.

\newtheorem*{nonsepThm}{Theorem \ref{nonsep}}
\newcommand\nonsep{ Let $M$ be a closed, simple 3-manifold
  containing a nonseparating connected closed incompressible surface
  $S$ of genus $g$.  Suppose that $\bar{\chi}(\kish(M \cut S)) <
  2g-2$, and that $M$ is $(2g-1,g)$-small. Then
  $\Hg(M)\le2g+1$. }
\begin{nonsepThm}
\nonsep
\end{nonsepThm}

In a closed, simple $3$-manifold, every connected closed
incompressible surface has genus at least $2$. Thus any such manifold
is $(2,2)$-small.  Hence applying Theorems \ref{trichotomy} and
\ref{nonsep} to a manifold containing an embedded surface of genus
$2$, we will easily obtain the following corollary.

\newtheorem*{genustwoCor}{Corollary \ref{genustwo}}
\newcommand\genustwo{ Suppose that $M$ is a closed, simple
  $3$-manifold which contains a connected closed incompressible
  surface of genus $2$, and that $\Hg(M)\ge6$. Then $M$ contains a
  connected closed incompressible surface $S$ of genus $2$ such that
  either $M \cut S$ has an acylindrical component, or
  $\chibar(\kish(M\cut S))\ge2$.  }
\begin{genustwoCor}
\genustwo
\end{genustwoCor}

This corollary will suffice for the geometric applications in this
paper.  In a future paper, we will apply Theorems \ref{nonsep} and
\ref{trichotomy} to the case of a genus 3 surface.  

In combination
with work of Agol-Storm-Thurston \cite{ASTD} and Kojima-Miyamoto
\cite{KM}, Corollary \ref{genustwo} implies the following volume bound
for sufficiently complex hyperbolic Haken manifolds.

\newtheorem*{manganeseThm}{Theorem \ref{manganese}}
\newcommand\manganese{
  Let $M$ be a closed, orientable hyperbolic 3--manifold containing a
  closed connected incompressible surface of genus $2$, and
  suppose that $\Hg(M)\ge6$.  Then $M$ has volume greater than $6.45$.
}
\begin{manganeseThm}
\manganese
\end{manganeseThm}

Theorem \ref{manganese} implies Theorems \ref{volume5.06} 
and \ref{volume3.08} below, which extend earlier work of Culler-Shalen.

\newtheorem*{fiveosixThm}{Theorem \ref{volume5.06}}
\newcommand\fiveosix{ 
  Let $M$ be a complete, finite--volume, orientable hyperbolic
  3-manifold having exactly one cusp, and suppose that
  $$\mathrm{dim}_{\mathbb{Z}_2} H_1(M;\mathbb{Z}_{2}) \geq 6.$$ 
 Then $M$ has volume greater than $5.06$.
}
\begin{fiveosixThm}
\fiveosix
\end{fiveosixThm}

Theorem \ref{manganese} is an improvement on \cite[Proposition
10.1]{CS_onecusp}.  There, the stronger lower bound of $7$ on the
dimension of $\mathbb{Z}_{2}$--homology gives only a conditional
conclusion: either the volume bound above holds, or $M$ contains an
embedded connected closed incompressible surface of genus $2$. The
weakening of the lower bound on the dimension of homology in the
hypothesis is made possible by the results of \cite{CS_vol}. In the
case where $M$ contains an embedded connected closed incompressible
surface of genus $2$ and $H_1(M;\mathbb{Z}_2)$ has dimension at least
$6$, a Dehn filling argument combined with Theorem \ref{manganese}
gives a better volume bound of $6.45$.
  
\newtheorem*{threeoeightThm}{Theorem \ref{volume3.08}}
\newcommand\threeoeight{Let $M$ be a closed, orientable hyperbolic
  3--manifold with
  $$\mathrm{dim}_{\mathbb{Z}_2} H_1(M; \mathbb{Z}_2) \geq 4$$ 
  and suppose that the image of the cup product map  $H^1 (M;
  \mathbb{Z}_{2})\otimes H^1 (M;
  \mathbb{Z}_{2})\to H^2 (M;
  \mathbb{Z}_2)$ has dimension at most $1$.  Then $M$ has volume
  greater than $3.08$.
}
\begin{threeoeightThm}
\threeoeight
\end{threeoeightThm}

Theorem \ref{volume3.08} should be compared with Theorem \geomsix\ of
\cite{CS_vol}, which gives the same conclusion under the hypothesis
that the $\mathbb{Z}_2$--homology of $M$ has dimension at least $6$,
and with no restriction on cup product. As with that theorem, the
proof of Theorem \ref{volume3.08} uses the fact that if $\pi_1(M)$ is
$3$--free, $M$ has volume greater than $3.08$ (see Corollary 9.3 of
\cite{last}).  If $\pi_1(M)$ has a $3$--generator subgroup $G$ which
is not free, the homological hypotheses of Theorem \ref{volume3.08}
ensure that $M$ has a two-sheeted cover $\widetilde{M}$ to which $G$
lifts, with $\mathrm{dim}_{\mathbb{Z}_2} H_1(M;\mathbb{Z}_2) \geq 6$.
Then Theorem 1.1 of \cite{CS_vol} implies that $\widetilde{M}$
contains a connected closed incompressible surface of genus
2, and Theorem \ref{manganese} implies that
$\widetilde{M}$ has volume greater than $6.45$, hence that $M$ has
volume greater than $3.22$.

\section{Topological preliminaries}\label{sec:prelim}

In general we will follow \cite{hempel} for standard terminology
concerning $3$-manifolds. (This includes, for example, the terms
``irreducible'' and ``boundary-irreducible'' which were used in the 
introduction.)  Here we will explain a few
special conventions and collect some preliminary results used 
throughout this paper.

We will work in the PL category in Sections
\ref{sec:prelim}---\ref{sec:trichotomy}, and in the smooth category in
Section \ref{sec:volumes}. The only result from the earlier sections
quoted in Section \ref{sec:volumes} is Corollary \ref{genustwo}, and
the smooth version of this result follows from the PL version. We will
also use, generally with explicit mention, the well-known fact that a
closed, orientable hyperbolic $3$-manifold is simple.

In Sections \ref{sec:prelim}---\ref{sec:trichotomy} we will use the
following conventions concerning regular neighborhoods. Let $K$ be a
compact polyhedron in a PL $n$-manifold $M$. We define a {\it
  semi-regular neighborhood} of $K$ in $M$ to be a neighborhood of $K$
which is a compact PL submanifold of $M$ and admits a polyhedral
collapse to $K$. We define a {\it regular neighborhood} of $K$ in $M$
to be a semi-regular neighborhood $N$ of $K$ in $M$ such that
$N\cap\partial M$ is a semi-regular neighborhood of $K\cap\partial M$
in $\partial M$.

Let $Y$ be a subset of a topological space $X$, and suppose that $X$
and $Y$ are locally path connected. We will say that $Y$ is {\it $
  \pi_1$-injective} in $X$ if whenever $A$ and $B$ are components of
$X$ and $Y$ respectively, such that $B\subset A$, the inclusion
homomorphism $\pi_1(B)\rightarrow \pi_1(A)$ is injective.

A closed orientable surface $S$ in the interior of an orientable
$3$-manifold $M$ will be termed {\it incompressible} if $S$ is
$\pi_1$-injective in $M$ and no component of $S$ is a sphere. We shall
not use the term ``incompressible'' for bounded surfaces.

We follow the conventions of \cite{ScharThom} regarding
Heegaard splittings and compression bodies.  The following standard
fact is a direct consequence of the definitions.

\Lemma\label{Heegaard genus is easy}  Let $Q$ be an orientable 
$3$-manifold with boundary, and suppose $\cals$ is a 
Heegaard surface in $Q$.  \begin{enumerate}
\item  Let $Q'$ be obtained from $Q$ by attaching a 
$2$-handle to a component of $\partial Q$. Then $\cals$ is a
Heegaard surface in $Q'$.
\item  Let  $Q'$ be obtained from $Q$ by attaching a handlebody 
to $Q$ along a  component of $\partial Q$.  Then $\cals$ is a
Heegaard surface in $Q'$.
\end{enumerate} \EndLemma
\NoProof

The lemma below is also standard, and will be used in 
Sections \ref{sec:nonsep} and \ref{sec:trichotomy}.

\Lemma \label{Heegaard genus is haard}  Let $Q$ be an
orientable $3$-manifold with boundary, and suppose $\cals$ 
is a Heegaard surface in $Q$ of genus $g$.  Let
$Q'$ be obtained from $Q$ by adding a $1$-handle with both attaching
disks in the same component of $Q-\cals$. Then $\Hg(Q')\le  g+1$.
\EndLemma

\Proof 
By definition we have $Q = C_1\cup C_2$ where $C_1$ and $C_2$ are
compression bodies such that $\partial_+ C_1 = \cals = \partial_+ C_2$
and $C_1\cap C_2 = \cals$.  After relabeling we may assume that $Q'$
is obtained by attaching a $1$-handle $H$ to $\partial_- C_1$.  We may
write $C_1$ as $(\cals \times I) \cup \mathcal{D} \cup \mathcal{T}$
where $\cals$ is identified with $\cals\times\{1\}$, $\mathcal{D}$ is
a union of disjoint $2$-handles attached along annuli in $\cals \times
\{0\}$, and $\mathcal{T}$ is a union of $3$-handles.  Since
$\mathcal{D} \cap \partial C_1$ is a union of disjoint disks, there is
an ambient isotopy of $C_1$ which is constant on $\cals$ and which
moves the two attaching disks of $H$ so that they are disjoint from
$\mathcal{D}$.  We may thus assume that the attaching disks for $H$
are contained in $\cals\times\{0\}$.

Let $N_0$ be a regular neighborhood in $H$ of its
core.  We have $N_0 \cap \partial_- C_1 = (E \times \{0\}) \cup (E'
\times \{0\})$, where $E$ and $E'$ are disjoint disks in $\cals =
\partial_+ C_1$.  Let
 $N = N_0 \cup (E \times I) \cup (E' \times I)$,
so that $E = E\times\{1\}$ and $E' = E'\times\{1\}$ are contained
in $\partial N$.
Set
$$\cals' = \overline{\cals-(E \cup E')} \cup
           \overline{\partial N - (E \cup E')}.$$ 
The surface $\cals'$ has genus $g+1$ by construction.  To complete the
proof we will show that it is a Heegaard surface for $Q \cup H$.

Let $P = (\overline{\cals - (E \cup E')} \times I)
\cup \overline{H - N_0}$ and set $C_1' = P \cup \mathcal{D}$.
Note that $P$ is a semi-regular neighborhood of $\cals'$ in $C_1'$,
and hence homeomorphic to $\cals' \times I$.  The attaching
annuli of the $2$-handles in $\mathcal{D}$ lie in the component of
$\partial P$ which is disjoint from $\cals'$.  It follows that
$C_1'$ is a compression body with $\partial_+C_1' = \cals'$.

Next let $C_2' = C_2 \cup N$.  From the dual description of $C_2$ as
$(\partial_- C_2)\times I$ with a collection of
$1$-handles attached, it follows that $C_2 \cup N$ is a compression
body with $\partial_+C_2' = \cals'$.

By construction we have $Q' = C_1'\cup C_2'$ and $C_1'\cap C_2'=
\partial_+C_1' = \partial_+C_2' = \cals'$.  Thus $\cals'$ is indeed
a Heegaard surface for $Q'$.
\EndProof

The following relatively straightforward result will be used in
Sections \ref{sec:nonsep} and \ref{sec:trichotomy}.

\Proposition\label{misty}
Let $g\ge2$ be an integer. Let $M$ be an irreducible, orientable
$3$-manifold which is $(g,\frac{g+1}{2})$-small.  Let $V$ be
a compact, connected, irreducible $3$-dimensional submanifold of $M$
which is $\pi_1$-injective. Suppose that either \Alternatives
\item\label{won't look like rain} $\chibar(V)<g-1$, or
\item\label{won't look like snow} $\chibar(V)\le g-1$ and $V$ is 
boundary-reducible. \EndAlternatives
Then $V$ is a handlebody.
\EndProposition

\Proof
Choose a properly embedded (possibly empty) submanifold $\cald$ of $V$
such that
\begin{enumerate}
\item\label{our magic can do anything} each component of $\cald$ is an
essential disk,
\item\label{here and now boys} no two components of $\cald$ are parallel, and
\item $\cald$ is maximal among all properly embedded submanifolds of
$V$ satisfying (\ref{our magic can do anything}) and (\ref{here and
now boys}).
\end{enumerate}
(Since $V$ is irreducible, (\ref{our magic can do anything}) and
(\ref{here and now boys}) imply that the components of
$\partial\cald$ are non-trivial and pairwise non-parallel simple
closed curves in $\partial V$; hence a submanifold $\cald$ satisfying
(\ref{our magic can do anything}) and (\ref{here and now boys}) has at
most $3\chibar(V)$ components, and hence a maximal submanifold with
these properties exists.)

Let $N$ be a regular neighborhood of $\cald$ in $V$, and set
$Q=\overline{V-N}$. In order to complete the proof it suffices to show
that every component of $Q$ is a ball.

Let us denote by $B_1,\ldots,B_\nu$ the components of $Q$ that are
balls, and by $R_1,\ldots,R_k$ the remaining components of $Q$. A
priori we have $k \ge0$ and $\nu\ge0$. We must show that $k=0$.

If $Q$ contains an essential disk $D$ we may assume after an isotopy
that $D\cap N=\emptyset$; then $D\cup\cald$ is a properly embedded
submanifold of $V$ satisfying (\ref{our magic can do anything}) and
(\ref{here and now boys}), a contradiction to the maximality of
$\cald$. This shows that $Q$ is boundary-irreducible, so that
$\partial Q$ is $\pi_1$-injective in $Q$. But $Q$ is $\pi_1$-injective
in $V$ since every component of $\cald$ is a disk, and $V$ is
$\pi_1$-injective in $M$ by hypothesis. Hence: \Claim\label{i slept}
$\partial Q$ is $\pi_1$-injective in $M$.
\EndClaim

The manifold $Q$ is obtained from the irreducible
$V$ by splitting along a collection of disjoint properly embedded
disks. Hence:

\Claim\label{in an oxygen tent} Each component of $Q$ is  irreducible.
\EndClaim

It follows from \ref{in an oxygen tent} that   each $R_i$ is
irreducible. Since  by definition no $R_i$ is  a ball, we deduce:

\Claim\label{life for people who can't read}
No boundary component of any $R_i$ can be a sphere.
\EndClaim

Let $n$ denote the number of components of $\cald$, and observe that
$\chibar(Q)=\chibar(V)-n$. Next we note that by properties (\ref{our
  magic can do anything}) and (\ref{here and now boys}) of $\cald$,
each component of $ Q$ which is a ball must contain at least three
components of the frontier $F$ of $N$ in $V$. Since each component of
$N$ contains exactly two components of $F$, the number $\nu$ of
components of $Q$ that are balls is at most $2n/3$. Hence
$$\begin{aligned}\chibar(V)-n&=\chibar(Q)\\
&=-\nu+\sum_{i+1}^k\chibar(R_i)\\
&\ge-(2n/3)+\sum_{i+1}^k\chibar(R_i),
\end{aligned}$$
so that 
\Equation\label{for the fbi}
\sum_{i+1}^k\chibar(R_i)\le\chibar(V)-(n/3).
\EndEquation
If alternative (\ref{won't look like rain}) of the hypothesis
holds we have $\chibar(V)<g-1$ and $n\ge0$; if alternative
(\ref{won't look like snow}) holds, we have $\chibar(V)\le g-1$ and
$n>0$. Thus in either case, (\ref{for the fbi}) implies that
\Equation\label{and found god}
\sum_{i+1}^k\chibar(R_i)<g-1.
\EndEquation

On the other hand, \ref{life for people who can't read} implies that
$\chibar(R_i)\ge0$ for each $i$ with $1\le i\le k$. In view of
(\ref{and found god}) it follows that
\Equation\label{and i mean it}
\chibar(R_i)< g - 1
\EndEquation
for each $i$ with $1\le i\le k$. 

We now proceed to the proof that $k=0$. Suppose that $k\ge1$, and
consider the manifold $R_1$. By \ref{i slept}, the boundary of $R_1$
is $\pi_1$-injective in $M$. By \ref{life for people who can't read},
no component of $\partial R_1$ is a sphere. Hence every component of
$\partial R_1$ is incompressible in $M$. Furthermore, we have
$$\chibar(\partial R_1)=2\chibar(R_1)<2g-2$$
by (\ref{and i mean it}). If $\partial R_1$ is connected, it is a
separating connected closed incompressible surface with
$\chibar(\partial R_1)<2g-2$, so that its genus is strictly less than
$g$. This contradicts the hypothesis.

Now suppose that $\partial R_1$ is disconnected. Since
$\chibar(\partial R_1)<2g-2$ and no component of $\partial R_1$ is a
sphere, we have $\chibar(S)<g-1$ for some component $S$ of $\partial
R_1$. This means that $S$ is a connected closed incompressible surface
of genus strictly less than $\frac{g+1}{2}$. Again we have a
contradiction to the hypothesis.
\EndProof

\section{Non-separating surfaces}\label{sec:nonsep}

The purpose of this section is to prove Theorem \ref{nonsep}, which
was stated in the introduction.

\Definition \label{horvert}  If a $3$-manifold $X$ has the
structure of an $I$-bundle over a surface $T$ and $p : X \rightarrow T$ is
the bundle projection, we will call $\partial_v X \doteq p^{-1}(\partial T)$ 
the \textit{vertical} boundary of $X$ and 
$\partial_h X \doteq \overline{\partial X - \partial_v X}$ the 
\textit{horizontal} boundary of $X$.  \EndDefinition 

Note that $\partial_v X$ inherits 
the structure of an $I$--bundle over $\partial T$, and $\partial_h X$ 
 the structure of a $\partial I$-bundle over $T$, from the original
$I$-bundle structure on $X$.

\Theorem \label{nonsep}
\nonsep
\EndTheorem

\Proof
Let $M$ and $S$ be as in the statement of the theorem. Set
$M'=M \cut S$, and note that $\bar{\chi}(M') = 2g-2$.  Since by
hypothesis we have $\bar{\chi}(\kish(M')) < 2g-2$, the characteristic
submanifold of $M'$ has a component $X$ which is an $I$--bundle
over a surface with negative Euler characteristic.  We
identify $M'$ with $\overline{M-N}$, where $N$ is a regular
neighborhood of $S$ in $M$; we then have
$X\subset\overline{M-N}\subset M$. We set $\Sigma=N\cup X\subset
M$.  Since the horizontal boundary of $X$ has Euler characteristic
$2\chi(X)$, we have
$$ \chibar(\Sigma) = \chibar(N) + \chibar (X) - 2\chibar(X)
= 2g - 2- \chibar(X).$$ 
Since $\chi(X)<0$, it follows that
$\chibar(\Sigma)< 2g-2$. 

Set $K =\overline{M-\Sigma}\subset\overline{ M-N}=M'$. (It may
happen that $K=\emptyset$.) Since $\partial
K=\partial\Sigma$, we have
$\chibar(K)=\chibar(\Sigma)$, and hence
\Equation\label{you won't catch me}
\chibar(K)< 2g-2.
\EndEquation

Since the frontier components of $K$ in $M'$ are essential
annuli, $K$ is $\pi_1$-injective in $M'$. The
incompressibility of $S$ implies that $M'=\overline{M-N}$ is
$\pi_1$-injective in $M$. Hence:
\Claim\label{old slithergadee}
$K$ is $\pi_1$-injective in $M$. 
\EndClaim
Note also that $M'$ is irreducible because the surface $S$ is
incompressible in the irreducible $3$-manifold $M$. The manifold $K$
is a union of components of the manifold obtained by splitting $M'$
along a collection of disjoint properly embedded annuli. Hence:
\Claim\label{you may catch all the others}
Every component of $K$ is irreducible.
\EndClaim
Since each component of $K$ contains a 
component of the frontier of the characteristic submanifold of $M'$,
which is an essential annulus in $M'$, no component of $K$ is a
ball. In view of \ref{you may catch all the others} it follows that no
component of $\partial K$ is a sphere. Hence:
\Claim\label{but you wo---}
Every component
of $K$ has  non-positive Euler characteristic. 
\EndClaim
Now consider any component $V$ of $K$. Set $V'=K-V$. It follows from
\ref{but you wo---} that $\chibar(V')\ge0$.  We have
$\chibar(V)=\chibar(K)-\chibar(V')$, and hence by
\ref{you won't catch me} 
\Equation\label{savonarola brown}
\chibar(V)< 2g-2. 
\EndEquation
By hypothesis  $M$ is $(2g-1,g)$-small. Since
$g=\frac{(2g-1)+1}{2}$, this means that $M$ is
$(2g-1,\frac{(2g-1)+1}{2})$-small. In view of 
\ref{old slithergadee},
\ref{you may catch all the others} and \ref{savonarola brown}, case
(\ref{won't look like rain}) of the hypothesis of Proposition
\ref{misty} holds with $2g-1$ playing the role of $g$. Proposition
\ref{misty} therefore implies that $V$ is a handlebody.

Thus we have shown:
\Claim\label{i hate to say i told you so}Every component of 
$K =\overline{M-\Sigma}$  is a handlebody.
\EndClaim

We now turn to the estimation of $\Hg(M)$. First note that since $N$
is a trivial $I$-bundle over a surface of genus $g$, it can be
obtained from a handlebody $J$ of genus $2g$ by adding a $2$-handle.
The boundary $\cals$ of a collar neighborhood of $\partial J$ in $J$
is a Heegaard surface of genus $2g$ in $J$. Hence by assertion (1) of
Lemma \ref{Heegaard genus is easy}, $\cals$ is a Heegaard surface in
$N$. Note that $\partial N$ is contained in a single component of
$N-\cals$.

On the other hand, recall that $\Sigma=N\cup X$, where $X$ is an
$I$-bundle over a connected surface $T$ and $N\cap X = \partial_h X$ 
is the horizontal boundary of $X$.  Let $E\subset T$ be a disk such that
for each boundary component $c$ of $T$, the set $E\cap c$ is a
non-empty union of disjoint arcs in $c$. Let $p:X\to T$ denote the
bundle projection, and set $Y=p^{-1}(E)$.  Then $Y$ inherits the 
structure of a (necessarily trivial) $I$-bundle over $E$, and 
$Y \cap N = Y \cap \partial_h X$ is the horizontal boundary of $Y$,
consisting of two disks.  Thus the set $Y$ may be thought of 
as a $1$-handle
attached to the submanifold $N$.  Since $\cals$ is a genus--$2g$
Heegaard surface in $N$, and $\partial N$ is contained in a single
component of $N-\cals$, it follows from Lemma \ref{Heegaard genus is
 haard} that $\Hg(N\cup Y)\le2g+1$.

Next, note that each component of $\overline{(\partial T)-E}$ is an
arc, and hence that each component of the set $\cald\doteq
p^{-1}(\overline{(\partial T)-E)})$ is a disk. Note also that
$\cald\cap(N\cup Y)=\partial\cald$. Hence if $R$ denotes a regular
neighborhood of $\cald$ relative to $\overline{X-Y}$, saturated in the
fibration of $X$, the manifold $N\cup Y\cup R$ is obtained from $N\cup
Y$ by adding finitely many $2$-handles. By assertion (1) of Lemma
\ref{Heegaard genus is easy} it follows that
$$\Hg(N\cup Y\cup R)\le\Hg(N\cup Y)\le2g+1.$$

Finally, note that each component of $\overline{M-(N\cup Y\cup R)}$ is
either
\begin{enumerate}[(a)]
 \item a component of $\overline{X-(Y\cup R)}$ or
 \item a component of $\overline{M-\Sigma}$.
\end{enumerate}
Each component of type (a) is a sub-bundle of $X$ over a bounded
subsurface, and is therefore a handlebody.  Each component of type (b)
is a handlebody by virtue of \ref{i hate to say i told you so}. Since
each component of $\overline{M-(N\cup Y\cup R)}$ is a handlebody, it
now follows from assertion (2) of Lemma \ref{Heegaard genus is
easy} that
$$\Hg(M)\le\Hg(N\cup Y\cup R)\le2g+1.$$
\end{proof}

\section{Annulus bodies and \sosmall\ manifolds} \label{sec:annulus}

In the next two sections, we develop an organizing principle for
cylinders properly embedded in the complement of a separating
connected closed incompressible surface.

\Definition\label{pineapple}
Let $Y$ be a compact, connected $3$-manifold, and let $S$ be a
  (possibly disconnected) closed, $2$-dimensional submanifold of
  $\partial Y$. We shall say that $Y$ is
an \textit{annulus  body} relative to $S$ if there is a
  properly embedded annulus $A\subset Y$ with $\partial A\subset S$,
  such that $Y$ is a semi-regular neighborhood of $S\cup A$.
\EndDefinition

\Lemma\label{re-verse it}Let $Y$ be a compact, connected
  $3$-manifold, and let $S\subset\partial Y$ be a closed
  $2$-manifold. If $Y$ is
an annulus  body relative to $S$, then $Y$ is also
an annulus  body relative to $(\partial Y)-S$. Furthermore,
we have $\chibar(Y)=\chibar(S)=\chibar((\partial Y)-S)$.
\EndLemma

\Proof
We set $T=(\partial Y)-S$.

By the definition of an annulus  body, $Y$ is a
semi-regular neighborhood of $S\cup A$ for some properly embedded
annulus $A\subset Y$ with $\partial A\subset S$.  Let $R$ be a regular
neighborhood of $A$ in $Y$.  Then there is a PL homeomorphism
$j:S^1\times [-1,1]\times [-1,1]\to R$ such that
$j(S^1\times\{0\}\times[-1,1])=A$ and
$j(S^1\times[-1,1]\times\{-1,1\})=R\cap S$. Let $B$ denote the annulus
$j(S^1\times[-1,1]\times\{0\})$. Set $Q=\overline{Y-R}$, and let $N$
be a regular  neighborhood of $Q\cap
S$ in $Q$, chosen small enough so that $N\cap B=\emptyset$. Set
$Y'=N\cup R$. Then $Y'$ is a compact $3$-manifold and
$S\subset\partial Y'$. If we set $T'=(\partial Y')-S$, then the
annulus $B$ is properly embedded in $Y'$ and $\partial B\subset
T'$. Furthermore, $Y'$ is a semi-regular neighborhood of $T'\cup B$, and
hence $Y'$ is an annulus  body relative to $T'$.

On the other hand, $Y$ and $Y'\subset Y$ are both semi-regular
neighborhoods of $S\cup A$, and $Y'\cap\partial Y=S$. Hence
$\overline{Y-Y'}$ is a collar neighborhood of $T\subset\partial Y$ in
$Y$. In particular the pairs $(Y,T)$ and $(Y',T')$ are PL
homeomorphic, and so $Y$ is an annulus  body relative to
$T$.

To prove the second assertion, we note that since $Y$ and $S\cup A$
are homotopy equivalent, we have $\chi(Y)=\chi(S\cup A)$; and that
since $A$ and $A\cap S=\partial A$ have Euler characteristic $0$, we
have $\chi(S\cup A)=\chi(S)$. This proves that
$\chibar(S)=\chibar(Y)$. Since we have shown that 
$Y$ is also an annulus  body relative to
$T$, we may substitute $T$ for $S$ in the last equality and conclude that
$\chibar(T)=\chibar(Y)$.
\EndProof

\Definition\label{leump}
Let $Z$ be a compact, connected, orientable $3$-manifold, and
let $S\subset\partial Z$ be a closed surface. We will say that $Z$
is \textit{\sosmall} relative to $S$ if $Z$ may be written in the form
$Z=Y\cup J$, where $Y\supset S$ and $J$ are compact
$3$-dimensional submanifolds of $Z$ such that
\begin{enumerate}
\item each component of $J$ is a handlebody,
	\item $Y$ is an annulus  body relative to
  $S$,
\item $Y\cap J= \partial J$, and
\item $\partial J$ is a union of components of $(\partial Y)-S$.
\end{enumerate}
(The submanifold $J$ may be empty.)
\end{definition}

\Lemma\label{take your leumps}Let $Q$ be a compact orientable
$3$-manifold, let $Z\subset Q$ be a compact submanifold whose frontier
$S$ is a connected closed surface in $\inter Q$, and suppose that $Z$ is \sosmall\ relative
to $S$. Then $\Hg(Q)\le1+\Hg(\overline{Q-Z})$.
\EndLemma

\Proof
We set $Q_0=\overline{Q-Z}$  and $g=\Hg(Q_0)$. 
We write $Z=Y\cup J$, where $Y$ and $J$ satisfy conditions
(1)--(4) of Definition \ref{leump}. Since
$Y$ is an annulus  body relative to
  $S$, it follows from Definition \ref{pineapple} that there is a
  properly embedded annulus $A\subset Y$ with $\partial A\subset S$,
  such that $Y$ is a semi-regular neighborhood, relative to $Y$
  itself, of $S\cup A$.

Let $\alpha$ denote a co-core of the annulus $A$, and fix a regular
neighborhood $h$ of $\alpha$ in $Y$ such that $h\cap A$ is
 a regular
neighborhood of $\alpha$ in $A$. The manifold $Q_0\cup h$ is obtained
from $Q_0$ by attaching a $1$-handle that has both its attaching disks
in the component $S$ of $\partial Q_0$. Hence it follows from
Lemma \ref{Heegaard genus is haard} that $\Hg(Q_0\cup h)\le1+\Hg(Q_0)=1+g$.

The disk $D=\overline{A-(h\cap A)}$ is properly embedded in the
manifold $\overline{Y- h}$. Hence if $R$ denotes a regular
neighborhood of $D$ relative to $\overline{Y- h}$, the manifold
$X_0\doteq Q_0\cup h\cup R$ is obtained from $Q_0\cup h$ by attaching
a $2$-handle. It therefore follows from assertion (1) of Lemma
\ref{Heegaard genus is easy} that $\Hg(X_0)\le\Hg(Q_0\cup h)\le1+g$.
But $X\doteq Q_0\cup Y$ is a semi-regular neighborhood of $X_0$
relative to $X$ itself, and is therefore homeomorphic to $X_0$. Hence
$\Hg(X)\le1+g$.

We have $Q=X\cup J$. In view of Conditions (1), (3) and (4) of
Definition \ref{leump}, it follows that each component of
$\overline{Q-X}$ is a handlebody whose boundary is contained in
$\partial X$. From assertion (2) of Lemma \ref{Heegaard genus
is easy} we deduce that $\Hg(Q)\le\Hg(X)\le1+g$.
\EndProof

\Lemma \label{small_homology}
Suppose that $Z$ is a compact, connected, orientable $3$-manifold,
that $\partial Z$ is connected, and that $Z$ is \sosmall\ relative to
$\partial Z$.  Let $g$ denote the genus of $\partial Z$. Then
$\Hg(Z)\le g+1$.
\EndLemma

\Proof
Let $N$ be a boundary collar for $Z$. Then $N$ has a Heegaard
splitting of genus $g$, the frontier $S$ of $N$ in $Z$ is connected,
and $\overline{Z-N}$ is \sosmall\ relative to $S$. The result
therefore follows upon applying Lemma
\ref{take your leumps}, with $Z$ and $\overline{Z-N}$ playing the
respective roles of $Q$ and $Z$ in that lemma.
\EndProof

\Lemma\label{i was a teenage werewolf}
Let $g\ge2$ be an integer. Let $Z$ be a compact, orientable
$3$-manifold having exactly two boundary components $S_0$ and $S_1$,
both of genus $g$. Then $Z$ is \sosmall\ relative to $S_0$ if and only
if either \Alternatives
\item\label{ranana banana} $Z$ is an annulus  body relative
to $S_0$, or
\item\label{surely bo-boorly}there is a solid torus $K\subset Z$ such
that $K\cap\partial Z$ is an annulus contained in $S_1$, and the pair
$(\overline{Z-K},S_0)$ is homeomorphic to $(S_0\times I,S_0\times\{0\})$.
\EndAlternatives
\EndLemma

\Proof
If alternative (\ref{ranana banana}) holds then $Z$ is
\sosmall\ relative to $S_0$: it suffices to take
$Y=Z$ and $J=\emptyset$ in Definition \ref{leump}. If alternative
(\ref{surely bo-boorly}) holds, let $J$ be a regular
neighborhood in $\inter K$ of a core curve of $K$ and set
$Y=\overline{Z-J}$. Then $Y$ is an annulus  body relative
to $S_1$. (The annulus $A$ appearing in Definition \ref{pineapple} is
bounded by two parallel simple closed curves.) It now follows from
Definition \ref{leump} that $Z$ is
\sosmall\ relative to $S_0$.

Conversely, suppose that $Z$ is \sosmall\ relative to $S_0$.  Let us
write $Z=Y\cup J$, where $Y$ and $J$ satisfy conditions (1)--(4) of
Definition \ref{leump}, with $S_0$ playing the role of $S$. Set
$T=(\partial Y)-S_0$.  Since $Y$ is an annulus body relative to $S_0$,
it follows from Lemma \ref{re-verse it} that $Y$ is an annulus body
relative to $T$. This means that $Y$ is a semi-regular neighborhood of
$T\cup A$, where $A$ is an annulus with $A\cap T=\partial A$. Since
$Y$ is connected it follows that $T$ has at most two components.

The conditions of Definition \ref{leump} imply that $T$ is the
disjoint union of $\partial J$ with $S_1$. Since $T$ has at most two
components and $S_1$ has exactly one component, $\partial J$ has at
most one component. If $\partial J=\emptyset$ then $J=\emptyset$,
i.e. $Z=Y$. This implies alternative (\ref{ranana banana}) of the
present lemma.

Now consider the case in which $\partial J$ has exactly one
component. In this case $J$ is a single handlebody, and $\partial J$
and $S_1$ are the components of $T$. According to Lemma \ref{re-verse
it} we have
$$2g-2=\chibar(S_0)=\chibar(T)=\chibar(S_1)+\chibar(\partial J).$$
But since $S_1$ is a surface of genus $g$ we have $\chibar(S_1)=2g-2$,
and hence $\chibar (J)=0$. Thus $J$ is a solid torus. 

Now $Y$ is a semi-regular neighborhood of $T\cup A=S_1\cup A\cup\partial
J$. Since $Y$ is connected, $A$ must have one boundary component in
$S_1$ and one in $\partial J$.  Let $R$ be a regular neighborhood
of $A$ in $Y$.  Then $R$ is a solid torus meeting $S_1$ and $\partial J$,
respectively, in regular neighborhoods of the simple closed curves
$A \cap S_1$ and $A \cap \partial J$.  Let $K = J \cup R$.  Since $J$ is a solid 
torus and $A \cap J$ is parallel in $R$ to its core, $K$ is a solid torus.
Furthermore, $K \cap S_1 = R \cap S_1$ is an annulus.

Now set $Q = \overline{Y- R} = \overline{Z-K}$.
If $N$ is a regular neighborhood in $Q$ of $Q \cap T$, then $Y ' = N \cup R$
is a semiregular neighborhood of $T \cup A$ contained in $Y$.  Therefore 
$\overline{Y - Y'}$ is a collar neighborhood of $S_0$ in $Y$; that is, the
pairs $(\overline{Y-Y'},S_0)$ and $(S_0 \times I,S_0 \times \{0\})$ are 
homeomorphic.  But by the definition of $N$, the pair 
$(\overline{Y - Y'},S_0) = (\overline{Q - N}, S_0)$
is homeomorphic to $(Q,S_0) = (\overline{Z-K},S_0)$, and
alternative (\ref{surely bo-boorly}) of the present lemma holds in this
case.
\EndProof

\Lemma\label{my father was a spy} Let $g\ge2$ be an integer. Let $Z$ be a compact, orientable
$3$-manifold having exactly two boundary components $S_0$ and $S_1$,
both of genus $g$. Then $Z$ is
\sosmall\ relative to $S_0$ if and
only if it is
\sosmall\ relative to $S_1$.
\EndLemma

\Proof By symmetry it suffices to show that if $Z$ is \sosmall\
relative to $S_0$ then it is \sosmall\ relative to $S_1$.  In view of
Lemma \ref{i was a teenage werewolf}, it suffices to show that if one
of the alternatives (\ref{ranana banana}) and (\ref{surely bo-boorly})
of that lemma holds, then it still holds when $S_0$ is replaced by
$S_1$. For alternative (\ref{ranana banana}) this follows from Lemma
\ref{re-verse it}. Now suppose that alternative (\ref{surely
  bo-boorly}) of Lemma \ref{i was a teenage werewolf} holds. Let $c$
be a core curve of the annulus $\overline{\partial K - (K\cap\partial
  Z)}$.  Since $(\overline{Z-K},S_0)$ is homeomorphic to $(S_0\times
I,S_0\times\{0\})$, there is a properly embedded annulus
$\alpha\subset \overline{Z-K}$ joining $c$ to a simple closed curve in
$S_0$. Let $B$ be a regular neighborhood of $K\cup\alpha$ in $Z$. Set
$P=\overline{Z-B}$, choose a regular neighborhood $N$ of $P\cup S_1$
in $Z$, and set $K'=\overline{Z-N}$. Then $K'$ is a solid torus,
$K'\cap\partial Z$ is an annulus contained in $S_0$, and the pair
$(\overline{Z-K'},S_1)=(N,S_1)$ is homeomorphic to $(S_1\times
I,S_1\times\{0\})$.
\EndProof

\section{Separating surfaces}  \label{sec:trichotomy}

\Number\label{i'll do a rhumba}
In this section we will use the theory of books of $I$-bundles as
developed in \cite{last}. We recall the definition here, in a slight
paraphrase of the form given in \cite{last}.

A {\it  book of $I$-bundles} is a triple
$\calw=(W,\calb,\calp)$, where $W$ is a (possibly disconnected) compact,
orientable $3$-manifold, and $\calb,\calp\subset W$ are submanifolds
such that
\Bullets
\item each component of $\calb$ is a solid torus;
\item $\calp$ is an $I$-bundle over a (possibly disconnected)
$2$-manifold, and every component of $\calp$ has Euler characteristic $\le0$;
\item $W=\calb\cup\calp$;
\item $\calb\cap\calp$ is the vertical boundary of $\calp$;
\item $\calb\cap\calp$ is $\pi_1$-injective in $\calb$; and
\item each component of $\calb$ meets at least one component of $\calp$.
\EndBullets

As in \cite{last}, we shall denote $W$, $\calb$ and $\calp$ by
$|\calw|$, $\calb_\calw $ and $\calp_\calw $ respectively. The
components of $\calb_\calw $ will be called {\it bindings} of $\calw$,
and the components of $\calp_\calw$ will be called its {\it pages}.
The submanifold $\calb\cap\calp$, whose components are properly
embedded annuli in $W$, will be denoted $\cala_\calw$.

An important observation, which follows from the definitions, is that
if $W$ is a simple  $3$-manifold with $\kish W=\emptyset$, then
$W=|\calw|$ for some book of $I$-bundles $\calw$.
\EndNumber

\Lemma\label{injective suboib}
If $\calw$ is any connected book of $I$-bundles then $\cala_\calw$ is
$\pi_1$-injective in $|\calw|$. Furthermore, $|\calw|$ is an
irreducible $3$-manifold.
\EndLemma

\Proof
If $A$ is any component of $\cala=\cala_\calw$, then $A$ lies
in the frontier of a unique component $P$ of $\calp=\calp_\calw$ and
in a unique component $B$ of $\calb=\calb_\calw$. Since $A$ is an
annulus of non-zero degree in the solid torus $B$, it is
$\pi_1$-injective in $B$. It is also $\pi_1$-injective in $P$, since
$A$ is a vertical boundary annulus of the $I$-bundle $P$ and
$\chi(P)\le0$. It follows that $\cala$ is $\pi_1$-injective in $\calb$
and in $\calp$ and hence in $W=|\calw|$, which is the first assertion. 

To prove the second assertion, we note that $\calb$ is irreducible
because its components are solid tori, and that $\calp$ is irreducible
because its components are $I$-bundles over surfaces of Euler
characteristic $\le0$. Thus $W$ contains the $\pi_1$-injective,
two-sided, properly embedded $2$-manifold $\cala$, and the manifold
obtained by splitting $W$ along $\cala$ is irreducible. It follows
that $W$ is itself irreducible.
\EndProof

\Lemma\label{small_Ibundles}
Let $M$ be a closed simple 3-manifold.  Suppose that $\calw$ is a
connected book of $I$-bundles with $W=|\calw|\subset M$, and that
$\partial W$ is a connected incompressible surface in $M$. Let $g$
denote the genus of $\partial W$. Suppose that $M$ is
$(g,\frac{g+1}{2})$-small.  Then $W$ is \sosmall\ relative to
$\partial W$.
\EndLemma

\Proof
We first consider the degenerate case in which $\calw$ has no bindings, so
that $W$ is an $I$-bundle over a closed surface $T$ with
$\chi(T)=2-2g<0$. We choose an orientation-preserving simple closed
curve $C\subset T$ and let  $A$ denote the annulus $p^{-1}(C)$, where
$p$ is the bundle projection. If $Y$ denotes a regular neighborhood of
$(\partial W)\cup A$, then $\overline{W-Y}$ is homeomorphic to an
$I$-bundle over a bounded surface, and hence to a handlebody; hence
$W$ is \sosmall\ in this case.

Now assume that $\calw$ has at least one binding, so that every page
is an $I$-bundle over a bounded surface.  The sum of the Euler
characteristics of the pages of $\calw$ is equal to
$\chi(W)=1-g<0$. In particular, $\calw$ has a page $P_0$ with
$\chi(P_0)<0$. Then $P_0$ is an $I$-bundle over some base surface $T$;
we let $p:P_0\to T$ denote the bundle projection.  We choose a
component $C$ of $\partial T$ and set $A =p^{-1}(C) \subset \partial_v
P_0$.  
(See Definition \ref{horvert}.)
Since $\chi(T)=\chi(P_0)<0$, there is an arc $\alpha\subset T$
such that $\partial\alpha=\alpha\cap\partial T \subset C$, and
$\alpha$ is not parallel in $T$ to an arc in $C$.  Now $A$ is a
properly embedded annulus in $W$, and $D=p^{-1}(\alpha)$ is a properly
embedded disk in $P_0$. The boundary of $D$ consists of two vertical
arcs in $A\subset\partial_vP_0$, and of two properly embedded arcs in
$\partial_hA_0$, each of which projects to $\alpha$ under the bundle
projection. Since $\alpha$ is not parallel in $T$ to an arc in $C$,
each of the four arcs comprising $\partial D$ is essential in either
$\partial_hP_0$ or $A$. Hence $\partial D$ is homotopically
non-trivial in $A\cup\partial_hP_0$.

We set $X=W\cut A$.  We may identify $P_0$ with a submanifold of
$X$. Since $\partial D\subset A\cup\partial_hP_0$, the disk $D$ is
properly embedded in $X$.  
Each component of $\partial_v P_0$
is an annulus in $\cala_\calw$ by definition.  Since the frontier curves of 
$A\cup\partial_hP_0$ in $\partial X$ are
boundary components of such annuli, it follows from 
Lemma \ref{injective suboib} that they are
homotopically nontrivial in $X$.
Therefore $A \cup \partial_h P_0$ is
$\pi_1$-injective in $\partial X$. Since $\partial D$ is homotopically
non-trivial in $A\cup\partial_hP_0$, it is homotopically non-trivial
in $\partial X$. Hence the disk $D$ is essential in $X$.

Note that $\partial X$ may have either one or two components. In
either case we have $\chibar(X)=\chibar(W)=g-1$. 

The manifold $W$ is irreducible by Lemma \ref{injective suboib}, and 
$X$ is obtained by splitting $W$ along a properly embedded surface. Hence:

\Claim\label{if not who, then me}
Each component of $X$ is irreducible.
\EndClaim

We shall identify $X$ homeomorphically with $\overline{W-N}$, where
$N$ is a regular neighborhood of $A$ in $W$. With this identification,
it follows from Lemma \ref{injective suboib} that $X$ is
$\pi_1$-injective in $W$. On the other hand, the incompressibility of
$\partial W$ implies that $W$ is $\pi_1$-injective in $M$. Hence:

\Claim\label{won't look like fog}
$X$ is $\pi_1$-injective in $M$.
\EndClaim

We now claim:
\Claim\label{that's all we know}
Each component of $X$ is a handlebody.
\EndClaim

To prove \ref{that's all we know}, we recall that by hypothesis $M$ is
$(g,\frac{g+1}{2})$-small.  Furthermore each component of $X$
is $\pi_1$-injective in $M$ by \ref{won't look like fog}, and is
irreducible by \ref{if not who, then me}. Hence it suffices to show
that one of the conditions (\ref{won't look like rain}) or (\ref{won't
  look like snow}) of Proposition \ref{misty} holds for each component
$V$ of $X$.

If $X$ is connected, then it is boundary-reducible since it contains
the essential disk $D$. Since $\chibar(X)=g-1$, condition (\ref{won't
  look like snow}) of Proposition \ref{misty} holds with $V=X$.

Now suppose that $X$ has two components $X_0$ and $X_1$. Each of these
is a union of pages and bindings of $\calw$, and we may suppose them
to be indexed so that $P_0\subset X_0$. For $i=0,1$ we have
$\chibar(X_i)=\sum\chibar(P)$, where $P$ ranges over the pages
contained in $X_i$. By the definition of a book of $I$-bundles, each
term $\chibar(P)$ is non-negative. In particular we have
$\chibar(X_1)\ge0$, and since $\chibar(P_0)>0$ we have
$\chibar(X_0)>0$. On the other hand, we have
$$\chibar(X_0)+\chibar(X_1)=\chibar(X)=g-1.$$ It follows that
$\chibar(X_1)<g-1$, so that condition (\ref{won't look like rain}) of
Proposition \ref{misty} holds with $V=X_1$. On the
other hand, we have $\chibar(X_0)\le g-1$, and $X_0$ is
boundary-reducible since it contains the essential disk $D$. Thus condition (\ref{won't look like rain}) of
Proposition \ref{misty} holds with $V=X_0$. This proves \ref{that's
all we know}.

If $Y$ denotes a regular neighborhood of $(\partial W)\cup A$, then
$Y$ is an annulus  body relative to $\partial W$, and
$\overline{W-Y}$ is homeomorphic to $X$ and is therefore a
handlebody. By Definition
\ref{leump} it follows that $W$ is \sosmall.
\EndProof

\Definition\label{dopey def} 
Let $M$ be a closed orientable $3$-manifold, and let $g$ be an integer
$\ge2$. We define a {\it $g$-\pile\ } to be a finite sequence
$(Z_0,S_1,Z_1,\ldots,S_n,Z_n)$, where
\begin{itemize}
  \item $n$ is a strictly positive integer,
  \item $S_1,\ldots,S_n$ are disjoint separating incompressible
    surfaces in $M$ with genus $g$,
  \item $Z_0, Z_1, \ldots, Z_n$ are the closures of the components of
    $M - (S_1 \cup \hdots \cup S_n)$,
  \item $\partial Z_0 = S_1$, $\partial Z_n = S_n$, and $\partial Z_i
    = S_i \sqcup S_{i+1}$ for $0 < i < n$, and
  \item for $0<i<n$, $Z_i$ is \sosmall\ relative to $S_i$ and is not
    homeomorphic to $S_i \times I$.
\end{itemize}

We shall call the integer $n$ the {\it depth} of the $g$-\pile\
$(Z_0,S_1,Z_1,\ldots,S_n,Z_n)$.  We will say that a $g$-\pile\
$(Z_0',S_1',Z_1',\ldots,S'_{n'},Z'_{n'})$ is a {\it
(strict) refinement} of $(Z_0,S_1,Z_1,\ldots,S_n,Z_n)$ if
$(S_1,\ldots,S_n)$ is a (proper) subsequence of the finite sequence
$(S'_1,\ldots,S'_{n'})$. A $g$-\pile\ will be called {\it maximal} if
it has no strict refinement.
\EndDefinition

\Lemma\label{but all the horsemen knew her} 
Let $M$ be a closed orientable $3$-manifold, and let $g$ be an integer
$\ge2$.  If $(Z_0,S_1,Z_1,\ldots,S_n,Z_n)$ is a $g$-\pile, then
$(Z_n,S_n,Z_{n-1},\ldots,S_1,Z_0)$ is also a $g$-\pile.  Furthermore,
if $(Z_0,S_1,Z_1,\ldots,S_n,Z_n)$ is maximal, then
$(Z_n,S_n,Z_{n-1},\ldots,S_1,Z_0)$ is also maximal.
\EndLemma

\Proof
If $(Z_0,S_1,Z_1,\ldots,S_n,Z_n)$ is a $g$-\pile, then 
for $0<i<n$, since $Z_i$ is \sosmall\ relative to $S_i$, it follows
from Lemma \ref{my father was a spy} that $Z_i$ is \sosmall\ relative
to $S_{i+1}$. This implies the first assertion. 

To prove the second assertion, suppose that
$(Z_n,S_n,Z_{n-1},\ldots,S_1,Z_0)$ is not maximal, so that it
has a strict refinement
$(Z'_{n'},S'_{n'},Z'_{n'-1},\ldots,S'_1,Z'_0)$. Then
$(Z'_0,S'_1,\ldots,Z'_{n-1},S'_n,Z'_n)$ is a 
$g$-\pile\ according to the first assertion, and
is a strict refinement of $(Z_0,S_1,\ldots,Z_{n-1},S_n,Z_n)$.
Hence $(Z_0,S_1,\ldots,Z_{n-1},S_n,Z_n)$ is not maximal.
\EndProof

\Lemma\label{she was only a jockey's daughter} 
Let $M$ be a closed orientable $3$-manifold, and let $g$ be an integer
$\ge2$.  Suppose that $M$ is $(g,\frac{g}{2}+1)$-small.  If
$(Z_0,S_1,Z_1,\ldots,S_n,Z_n)$ is a maximal $g$-\pile, then $Z_n$ is
either \sosmall\ relative to $S_n$ or acylindrical.
\EndLemma

\Proof
Let us suppose that $Z_n$ is not acylindrical. Let $A$ be an essential
annulus in $Z_n$.  We set $X=Z_n\cut A$. We shall identify $X$
homeomorphically with $\overline{Z_n-N}$, where $N$ is a regular
neighborhood of $A$ in $Z_n$.

If $Y$ denotes a regular neighborhood of $(S_n)\cup A$, then
$Y$ is an annulus  body relative to $S_n$, and
$\overline{Z_n-Y}$ is ambiently isotopic to $X$. Hence in order to
show that $Z_n$ is \sosmall, it suffices to show that each component of
$\overline{Z_n-Y}$ is a handlebody.

Note that $\partial X$ has at most two components, and hence that $X$
has at most two components. We have
\Equation\label{mystic men who eat boiled owls}
g-1=\chibar(Z_n)=
\chibar(X)=\sum_V\chibar(V),
\EndEquation
 where $V$ ranges over the components of $X$. The essentiality of $A$
implies that $\chibar(F)\ge0$ for each component $F$ of $\partial X$,
and hence that $\chibar(V)\ge0$ for each component $V$ of $\partial
X$. It therefore follows from (\ref{mystic men who eat boiled owls})
that $\chibar(V)\le g-1$ for each component $V$ of $X$.

Since $A$ is an essential annulus, $X$ is $\pi_1$-injective in
$Z_n$. On the other hand, the incompressibility of $S_n$ implies that
$Z_n$ is $\pi_1$-injective in $M$. Hence:

\Claim\label{we are men of groans and howls}
$X$ is $\pi_1$-injective in $M$.
\EndClaim

Since $M$ is irreducible, the incompressibility of $S_n$ implies that
$Z_n$ is irreducible.  Since $A$ is a properly embedded annulus in
$Z_n$, we deduce:

\Claim\label{tell us what you wish o king}
Each component of $X$ is irreducible.
\EndClaim

We now claim:
\Claim\label{shuffle duffle muzzle muff}
For each component $V$ of $X$, some component of $\partial V$ is
compressible in $M$.
\EndClaim

To prove \ref{shuffle duffle muzzle muff}, we first consider the case
in which $\partial V$ is disconnected. In this case, since
$\chibar(\partial V)=2\chibar(V)\le2g-2$, there is a component $F$ of
$\partial V$ with $\chibar(F) < g$; hence the genus of $F$ is
strictly less than $\frac{g}{2} + 1$.  Since $M$ is
$(g,\frac{g}{2}+1)$-small, the surface $F$ must be compressible.

We next consider the case in which $\partial V$ is connected and
$\chibar(V)<g-1$. In this case we have $\chibar(\partial V)<2g-2$, so
that $\partial V$ has genus strictly less than $g$. Furthermore,
$\partial V$ separates $M$. Since $M$ is
$(g,\frac{g}{2}+1)$-small, the surface $\partial V$ must be
compressible.

There remains the case in which $\partial V$ is connected and
$\chibar(V)=g-1$. In this case we set $S=\partial V$ and observe that
$ S$ is a separating surface of genus $g$. We shall assume that $S$ is
incompressible and derive a contradiction.  Since
$X'\doteq\overline{Z_n-Y}$ is ambiently isotopic to $X$, some
component $V'$ of $X'$ is ambiently isotopic to $V$, and so
$S'\doteq\partial V'$ is a separating connected closed
incompressible surface of genus $g$.

We distinguish two subcases, depending on whether (a) $V$ is the only
component of $X$, or (b) $X$ has a second component $U$.  In subcase
(a), the boundary components of $Y$ are $S_n$ and $S'$. Since $Y$ is
an annulus  body relative to $S_n$, it is in particular
\sosmall\ relative to $S_n$. Furthermore, $Y$ cannot be homeomorphic
to $S_n\times I$, because it contains the annulus $A$, which is
essential in $Z_n$---and hence in $Y$---and has its boundary contained
in $S_n$. It now follows from Definition \ref{dopey def} that
$(Z_0,S_1,Z_1,\ldots,S_n,Y,S',X')$ is a $g$-\pile. This
contradicts the maximality of $(Z_0,S_1,Z_1,\ldots,S_n,Z_n)$.

In subcase (b), it follows from (\ref{mystic men who eat boiled owls})
that $\chibar(U)=0$. Since $\partial X$ has at most two components,
$\partial U$ is a single torus. The simplicity of $M$ implies that
$\partial U$ is compressible in $M$, and since $U$ is
$\pi_1$-injective in $M$ by \ref{we are men of groans and howls},
$\partial U$ cannot be $\pi_1$-injective in $U$. As $U$ is irreducible
by \ref{tell us what you wish o king}, it now follows that $U$ is a
solid torus. Hence the component $U'$ of $X'$ which is ambiently
isotopic to $U$ is a solid torus. According to Definition \ref{leump},
this implies that $Z\doteq Y\cup U'$ is \sosmall.  The boundary
components of $Z$ are $S_n$ and $S'$.  The \sosmall\ manifold $Z$
cannot be homeomorphic to $S_n\times I$, because it contains the
annulus $A$, which is essential in $Z_n$---and hence in $Z$---and has
its boundary contained in $S_n$. It now follows from Definition
\ref{dopey def} that $(Z_0,S_1,Z_1,\ldots,S_n,Z,S',V')$ is a
$g$-\pile. This contradicts the maximality of
$(Z_0,S_1,Z_1,\ldots,S_n,Z_n)$,

This completes the proof of \ref{shuffle duffle muzzle muff}.

Next we claim:
\Claim\label{we just can't tell you any more}
Each component of $X$ is boundary-reducible.
\EndClaim

In fact, if a component $V$ of $X$ were boundary-irreducible, then
$\partial V$ would be $\pi_1$-injective in $V$. In view of \ref{we are
  men of groans and howls} it would follow that $\partial V$ is
$\pi_1$-injective in $M$. But this contradicts \ref{shuffle duffle
  muzzle muff}. Thus \ref{we just can't tell you any more} is
established.

We now turn to the proof that each component of $X$ is a handlebody,
which will complete the proof of the lemma.

Let $V$ be any component of $X$. We have observed that $\chibar(V)\le
g-1$. By \ref{we are men of groans and howls}, $V$ is
$\pi_1$-injective in $M$, by \ref{tell us what you wish o king} it is
irreducible, and by \ref{we just can't tell you any more} it is
boundary-reducible. Since the hypothesis implies in particular that
$M$ is $(g,\frac{g+1}{2})$-small, case (\ref{won't look like snow})
of the hypothesis of Proposition \ref{misty} holds. It therefore
follows from Proposition \ref{misty} that $V$ is a handlebody.
\EndProof

\Lemma\label{mouse with a glandular condition}
Let $M$ be a closed orientable $3$-manifold, and let $g$ be an integer
$\ge2$. Suppose that $M$ is $(g,\frac{g}{2} + 1)$-small.  If
$(Z_0,S_1,Z_1,\ldots,S_n,Z_n)$ is a maximal $g$-\pile, then $Z_{0}$ is
either \sosmall\ relative to $S_1$ or acylindrical.
\EndLemma

\Proof
This is an immediate consequence of Lemmas \ref{but all the horsemen
  knew her} and \ref{she was only a jockey's daughter}.
\EndProof

\Theorem\label{trichotomy}
\trichotomy
\EndTheorem

\Proof
It follows from the Haken finiteness theorem \cite[Lemma
13.2]{hempel} that the set of all depths of $g$-\pile s in $M$ is
bounded. In particular any $g$-\pile\ has a refinement which is a
maximal $g$-\pile.

By hypothesis $M$ contains some separating connected closed
incompressible surface $T$ of genus $g$. If $X$ and $Y$ denote the
closures of the components of $M-T$, then $(X,T,Y)$ is a $g$-\pile\ of
depth $1$. In particular $M$ contains a $g$-\pile, and hence contains
a maximal $g$-\pile.

Now suppose that the conclusion of Theorem \ref{trichotomy} does not
hold. Fix a maximal $g$-\pile\ $(Z_0,S_1,Z_1,\ldots,S_n,Z_n)$.  Then
neither $Z_0$ nor $Z_n$ is acylindrical, since otherwise $S = S_1$ or
$S_n$ would satisfy conclusion (1) of the theorem.  In view of
\ref{she was only a jockey's daughter} and \ref{mouse with a glandular
  condition}, and the hypothesis that $M$ is
$(g,\frac{g}{2}+1)$-small, it follows that $Z_0$ and $Z_n$
are both \sosmall.

For $0 < i \leq n$ we define $B_i^- = Z_0 \cup Z_1 \cup \hdots \cup
Z_{i-1}$ and $B_i^+ = Z_i \cup \hdots \cup Z_n$.  For each $i$, we
must have either $\kish(B_i^-)=\emptyset$ or $\kish(B_i^+)=\emptyset$,
since otherwise $S = S_i$ would satisfy conclusion (2) of the theorem.
Hence by the observation made in \ref{i'll do a rhumba}, at least one
of $B_i^-$ or $B_i^+$ has the form $|\calw|$ for some book of
$I$--bundles $\calw$.  Since $M$ is in particular
$(g,\frac{g+1}{2})$-small, it now follows from Lemma
\ref{small_Ibundles} that at least one of $B_i^-$ or $B_i^+$ is
\sosmall.

Since $B_0^-=Z_0$ is \sosmall, there is a largest index $k\le n$ such
that $B_k^-$ is \sosmall. Since $B_k^-$ is \sosmall, it follows from
Lemma \ref{small_homology} that $\Hg(B_k^-)\le g+1$. 
We distinguish two cases depending on whether $k<n$ or $k=n$.

If $k<n$ then $B_{k+1}^-$ is not \sosmall, and hence $B_{k+1}^+$ is
\sosmall.  By the definition of $B_{k+1}^-$, the frontier of 
$Z_k$ in $B_{k+1}^-$ is the closed surface $S_k$, and
by the definition of a $g$-\pile, $Z_k$ is \sosmall\ relative to $S_k$.
We may thus apply Lemma \ref{take your leumps} with $Q=B_{k+1}^-$
and $Z=Z_k$ to deduce that
$$\Hg(B_{k+1}^-)\le1+\Hg(B_k^-)\le g+2.$$

Now since
$B_{k+1}^+$ is \sosmall, we may again apply Lemma \ref{take your
leumps}, this time with $Q=M$ and $Z=B_{k+1}^+$, to deduce that
$$\Hg(M)\le1+\Hg(B_{k+1}^-)\le g+3.$$
This contradicts the hypothesis.

Now suppose that $k=n$. In this case, since $B_n^+=Z_n$ is \sosmall,
we may apply Lemma \ref{take your leumps} with $Q=M$ and $Z=Z_n$ to
deduce that
$$\Hg(M)\le1+\Hg(B_n^-)\le g+2.$$
Once again we have a contradiction to
the hypothesis.
\end{proof}

\Corollary\label{genustwo}
\genustwo
\EndCorollary

\Proof
Since $M$ is simple, it is $(2,2)$-small.

First consider the case in which $ M$ contains a separating,
connected, closed, incompressible surface of genus $2$. Since $M$ is
$(2,2)$-small, Theorem \ref{trichotomy} gives a separating
connected closed incompressible surface $S$ of genus $2$
such that either at least one component of $M\cut S$ is acylindrical,
or for each component $B$ of $M \cut S$ we have
$\kish(B)\ne\emptyset$. In particular, $\kish(M\cut S)$ has at least
two components. By Definition (see \ref{simple def}), each component of
$\kish(M\cut S)$ has a strictly negative Euler characteristic. Hence
$\chibar(\kish(M\cut S))\ge2$.

Now suppose that $ M$ contains no separating, connected, closed,
incompressible surface of genus $2$. In this case $M$ is
$(3,2)$-small. By hypothesis, $ M$ contains a connected, closed,
incompressible surface $S$ of genus $2$, which must be
non-separating. It now follows from Theorem \ref{nonsep} that 
$\bar{\chi}(\kish(M \cut S))\ge2$.
\EndProof

\section{Volume bounds} \label{sec:volumes}

Recall that a {\it slope} on a torus $T$ is an unoriented isotopy
class of homotopically non-trivial simple closed curves on $T$. If the
torus $T$ is a boundary component of an orientable $3$-manifold $N$,
and $r$ is a slope on $T$, we denote by $N(r)$ the ``Dehn-filled''
manifold obtained from the disjoint union of $N$ with $D^2\times S^1$
by gluing $(\partial D^2)\times S^1$ to $\partial N$ via a
homeomorphism which maps $(\partial D^2)\times\{{\rm point}\}$ to a
curve representing the slope $r$.

\Lemma\label{woo-woo ginsberg}
Let $N$ be a compact $3$-manifold whose boundary is a single torus,
let $S\subset N$ be a closed connected incompressible surface, and let
$p$ be a prime. Then there exist infinitely many slopes $r$ on
$\partial N$ for which the following conditions hold:
\begin{enumerate}
\item\label{herman had burped} the inclusion homomorphism $H_1(N;\ZZ_p)\to
  H_1(N(r_i);\ZZ_p)$ is an isomorphism; and
\item\label{the advantages of suburbia} $S$ is incompressible in $N(r_i)$.
\end{enumerate}
\EndLemma

\Proof
There is a natural bijective correspondence between slopes on
$\partial M$ and unordered pairs of the form $\{c,-c\}$ where $c$ is a
primitive element of $L\doteq H_1(\partial N;\ZZ)$. If $r$ is a
slope, the elements of the corresponding unordered pair are the
homology classes defined by the two orientations of a simple closed
curve representing $c$. If $c$ is a primitive element of $L$ we shall
denote by $r_c$ the slope corresponding to the pair $\{c,-c\}$.

Let $K\subset L$ denote the kernel of the natural homomorphism
$H_1(\partial N;\ZZ_p)\to H_1(N;\ZZ_p)$. If $c$ is a primitive class
$c\in K$, it follows from the Mayer-Vietoris theorem that the
inclusion homomorphism $H_1(N;\ZZ_p)\to H_1(N(r_c);\ZZ_p)$ is an
isomorphism.

We fix a basis $\{\lambda,\mu\}$ of $L$ such
that $\lambda\in K$, and we identify $L$ with an
additive subgroup of the two-dimensional real vector space
$V=H_1(\partial N;\RR)$. For each positive integer $n$, let 
$A_n\subset V$ denote the affine line $\RR\lambda+np\mu\subset
V$. Then $A_n\cap L\subset K$, and $A_n\cap L$ contains infinitely
many primitive elements of $L$ (for example the elements of the form
$(knp+1)\lambda+np\mu$ for $k\in\ZZ$). In particular, $K$ contains
infinitely many primitive elements of $L$.

We distinguish two cases, depending on whether there (a) does or (b)
does not exist an annulus in $M$ having one boundary component in $S$
and one in $\partial M$, and having interior disjoint from
$S\cup\partial M$. In case (a) it follows from \cite[Theorem
2.4.3]{cgls} that there is a slope $r_0$ such that for every slope $r$
whose geometric intersection number with $r$ is $>1$, the surface $S$
is incompressible in $M(r)$. (For the application of
\cite[Theorem 2.4.3]{cgls} we need to know that $S$ is not
boundary-parallel in $M$, but this is automatic since $S$ has genus
$2$.) In particular there are three affine lines $B_1$, $B_2$ and
$B_3$ in $V$ such that for any primitive class $c\in
L\setminus(B_1\cup B_2\cup B_3)$, the surface $S$ is incompressible in
$M(r_c)$.  If we choose a natural number $n$ large enough so that
$A_n$ is distinct from $B_1$, $B_2$ and $B_3$, then $A_n\cap(B_1\cup
B_2\cup B_3)$ consists of at most three points. Hence of the
infinitely many primitive elements of $L$ belonging to $A_n\cap L$, at
most three lie in $B_1\cup B_2\cup B_3$. For any primitive element
$c\in (A_n\cap L)\setminus(B_1\cup B_2\cup B_3)$, the slope $r_c$
satisfies Conclusions (\ref{herman had burped}) and (\ref{the
advantages of suburbia}).

In Case (b) it follows from \cite[Theorem 1]{wu} that there are at
most three slopes $r$ for which $S$
is compressible in $M(r)$. In particular, of the infinitely many
elements  $c\in K$ which are primitive in $L$, all but finitely many
have the property that $S$
is incompressible in $M(r_c)$. Hence there are infinitely many slopes
$c$ satisfying Conclusions (\ref{herman had burped}) and (\ref{the
advantages of suburbia}).
\EndProof

\Definition
If $X$ is a compact orientable manifold with non-empty boundary then by the
{\it double }of $X$ we shall mean the quotient space $DX$ obtained from
$X\times\{0,1\}$ by identifying $(x,0)$ with $(x,1)$ for each
$x\in\partial X$.  The involution of $X\times\{0,1\}$ which
interchanges $(x,0)$ and $(x,1)$ induces an orientation--reversing
involution $\tau:DX\to DX$, which we shall call the {\it canonical
involution} of $DX$.  Define $\mathrm{geodvol}\, X = \frac{1}{2}v_3\|[DX]\|$,
where $\|[DX]\|$ denotes the Gromov norm of the fundamental class
of $DX$, and $v_3$ is the volume of a regular ideal tetrahedron.
\EndDefinition

The following standard result does not seem to be in the literature.

\Proposition\label{totally geodesic}
Let $X$ be a compact connected orientable $3$-manifold with connected
boundary $S$ of genus greater than $1$.  Suppose that $X$ is irreducible,
boundary-irreducible and acylindrical.  Then $X$ admits a hyperbolic
metric with totally geodesic boundary, and $\mathrm{geodvol}\, X$ is the
volume of this metric.
\EndProposition

\Proof The closed manifold $DX$ is simple, and the surface $S$ is
incompressible in $DX$.  Thus $DX$ admits a complete hyperbolic metric
by Thurston's Hyperbolization Theorem for Haken manifolds \cite{otal}.
Let $\tau:DX\to DX$ be the canonical involution of $DX$.  Fix a
basepoint $\star\in S$ and a basepoint $\tilde\star$ in the universal
cover of $DX$ which maps to $\star$.  We identify $DX$ with
$\HH^3/\Gamma$, where $\Gamma$ is a Kleinian group.  Using the
basepoint $\tilde\star$, we identify $\pi_1(DX,\star)$ with $\Gamma$.
Let $p:\HH^3\to DX$ be the covering projection, and let $\tilde\tau
:\HH^3\to \HH^3$ denote the lift of $\tau\circ p$ which fixes
$\tilde\star$.  Let $\widetilde S$ be the component of $p^{-1}(S)$
which contains $\tilde\star$.  The map $\tilde\tau$ is then an
orientation-reversing involution of $\HH^3$ which fixes $\widetilde
S$.

Since $\tau$ is a homotopy equivalence, it follows from the proof of
Mostow's rigidity theorem (\cite{Mostow}, cf. \cite{BP}) that
$\tilde\tau$ extends continuously to $S^2_\infty$, and that there is
an isometry $\tau'$ of $\HH^3$ whose extension to $S^2_\infty$ agrees
with that of $\tilde\tau$.  In particular $\tau'$ is an
orientation-reversing isometry of $\HH^3$ whose restriction to
$S^2_\infty$ normalizes the restriction of $\Gamma$ to $S_2^\infty$.
Thus $\tau'$ normalizes $\Gamma$ in the isometry group of $\HH^3$, and
consequently $\tau'$ induces an involution of $DX$.  The restriction
of $\tau'$ to $S^2_\infty$ also commutes with the restriction of each
isometry in the image $\Delta$ of the inclusion homomorphism
$\pi_1(S,\star)\to\Gamma$.  It follows that $\tau'$ must be a
reflection through a hyperbolic plane $\Pi$, where $\Pi$ contains the
axis of each element of $\Delta$.  In particular $\Pi$ is invariant
under $\Delta$.  Moreover, since the image of $\Pi$ in $DX$ is
contained in the fixed set of an involution of $DX$, it must be a
compact subsurface $F$.

Since $F$ is covered by a hyperbolic plane, it is a totally geodesic
surface in $DX$.  Let $\Delta'\le\Gamma$ denote the stabilizer of
$\Pi$.  Then the covering space $\widetilde{DX} = \HH^3/\Delta'$ is
homeomorphic to $F\times\RR$.  Let $\widetilde F$ denote the image of
$\Pi$ in $\widetilde{DX}$, so that $\widetilde F$ is the image of a
lift of the inclusion of $F$ into $DX$.  Since $\Delta\le\Delta'$, the
inclusion of $S$ into $DX$ lifts to an embedding of $S$ in
$\widetilde{DX}$.  Let $\widetilde S$ denote the image of this lift.
The surfaces $\widetilde S$ and $\widetilde F$ are $\pi_1$-injective
and can be isotoped to a pair of disjoint surfaces which cobound a
compact submanifold $W$ of $\widetilde M$.  According to \cite[Theorem
10.5]{hempel} $W$ is homeomorphic to $\widetilde S\times I$.  This
shows that $\widetilde S$ is isotopic to $\widetilde F$ in
$\widetilde{DX}$, and hence that $S$ is homotopic to $F$ in $DX$.
Since $S$ is incompressible in $DX$, it follows from \cite[Corollary
5.5]{waldhausen} that $S$ is isotopic to the totally geodesic surface
$F$.

Since $S$ and $F$ are isotopic, there is an ambient isotopy of
$DX$ which carries $X$ onto a submanifold of $DX$ bounded by $F$.
Pulling back the hyperbolic metric by the time-$1$ map of this
isotopy endows $X$ with the structure of a complete hyperbolic
manifold with totally geodesic boundary.
The isometry of $DX$ induced by $\tau'$ fixes $F$ and exchanges its
complementary components.  Hence, 
with the pulled back metric, $X$ has
half the hyperbolic volume of $DX$.  Since the hyperbolic volume of
$DX$ is equal to $v_3\|[DX]\|$ (see \cite[Theorem C.4.2]{BP}), this
completes the proof.
\EndProof

The result below follows from a result of Agol-Storm-Thurston
\cite{ASTD}

\Proposition\label{ASTacylindrical}
Let $M$ be a closed, orientable hyperbolic 3-manifold containing a
closed connected incompressible surface $S$ such that $M
\cut S$ has an acylindrical component $X$.  Then
$$\mathrm{vol}\, M \geq \mathrm{geodvol}\, X.$$ 
\EndProposition

\Proof Let $g:S\to M$ denote the inclusion map.  Since $S$ is a
two-sided embedded surface, the family of all immersions of $S$ in $M$
which are homotopic to $g$ is non-empty.  Since $M$ is
$\mathbb{P}^2$-irreducible, the main result of \cite{FHS} asserts that
this family contains a least area immersion $f:S\to M$, which is
either an embedding or a $2$-sheeted covering of a non-orientable
surface.  Moreover, the second case arises only if $S$ bounds a
twisted $I$-bundle whose $0$-section is isotopic to $f(S)$.  Since $f$
is locally area minimizing, $f(S)$ is a minimal surface.

It follows from \cite[Corollary 5.5]{waldhausen} that if
$f$ is an embedding then $f(S)$ is ambiently isotopic to $S$, and if
$f$ is a $2$-sheeted covering map then $f(S)$ is isotopic to the
$0$-section of the twisted $I$-bundle bounded by $S$.  Hence one of
the components, say $X$, of $M\cut f(S)$ is acylindrical.
We may identify $X$ with the path completion of a component
$X_0$ of $M-f(S)$.  Then the natural map $X \to M$
maps the interior of $X$ homeomorphically onto $X_0$ and maps
$\partial X$ onto $f(S)$, either by a homeomorphism or by a
$2$-sheeted cover.  The latter possibility arises exactly when $S$
bounds a twisted $I$-bundle and $f(S)$ is a non-orientable
surface.  In particular, pulling back the hyperbolic metric on $M$
under the natural map $X\to M$ gives $X$ the structure of a complete
hyperbolic manifold whose boundary is a minimal surface.  Theorem 7.2
of \cite{ASTD} states that such a manifold $X$ satisfies $\vol X \ge
\mathrm{geodvol}\, X$.

Clearly we have $\vol M \ge \vol X$, so the proof is complete.
\EndProof

\Theorem\label{manganese}
\manganese
\EndTheorem

\Proof 
According to Corollary \ref{genustwo}, $M$ contains a connected closed
incompressible surface $S$ of genus $2$ such that either $M \cut S$
has an acylindrical component, or $\chibar(\kish(M\cut S))\ge2$.  If
$M \cut S$ has an acylindrical component $X$, by Proposition
\ref{totally geodesic} $X$ admits a hyperbolic metric with totally
geodesic boundary, and the volume of $X$ in this metric is equal to
$\mathrm{geodvol}\, X$.  The main result of \cite{KM} asserts that
this volume is greater than $6.45$; hence by Proposition
\ref{ASTacylindrical}, $\mathrm{vol}\, M$ is also greater than $6.45$.
On the other hand, if $\chibar(\kish(M\cut S))\ge2$, then Theorem 9.1
of \cite{ASTD} implies that $M$ has volume greater than $7.32$.
\EndProof

The following lemma is a strict improvement on Proposition 10.1 of
\cite{CS_onecusp}. The improvement is made possible by the results of
\cite{CS_vol}.

\Lemma\label{almost as much H_1}
Let $M$ be a complete, finite-volume, orientable hyperbolic $3$-manifold
having exactly one cusp, such that 
 $\dim_{\ZZ_2}H_1(M;\ZZ_2)\ge 6$. 
Then either
\begin{enumerate}
\item $\vol M>5.06$, or
\item $M$ contains a genus--$2$ connected incompressible surface.
\end{enumerate}
\EndLemma

\Proof
This is identical to the proof of
 \cite[Proposition
  10.1]{CS_onecusp} except that 
\Bullets
\item each of the two appearances of the
 number $7$ in the latter proof is replaced by $6$, and 
\item the
 reference to the case $g=2$ of \cite[Theorem \topeleven]{last} is
 replaced by a reference to the case $g=2$ of \cite[Theorem \topsix]{CS_vol}.
\EndBullets
\EndProof

\Theorem\label{volume5.06}
\fiveosix
\EndTheorem

\Proof For a hyperbolic manifold satisfying the hypotheses of Theorem
\ref{volume5.06}, Lemma \ref{almost as much H_1} asserts that either
$M$ has volume greater than $5.06$ or $M$ contains a closed connected
incompressible surface of genus $2$.  In the latter case, let $N$
denote the compact core of $M$.  According to Lemma \ref{woo-woo
  ginsberg} there is an infinite sequence of distinct slopes
$(r_i)_{i\ge1}$ on $\partial N$ such that $S$ is incompressible in
each $N(r_i)$, and $\dim_{\ZZ_2}H_1(N(r_i);\ZZ_2)\ge6$ for each $i$.
The hyperbolic Dehn surgery theorem (\cite{Th1}, cf. \cite{NZ})
asserts that $M_i\doteq N(r_i)$ is hyperbolic for all sufficiently
large $i$, and hence after passing to a subsequence we may assume that
all the $M_i$ are hyperbolic.  We now invoke Theorem 1A of \cite{NZ},
which implies that $\vol M_i < \vol M$ for all but finitely many $i$.
(The authors of \cite{NZ} attribute this particular consequence of
their main result to Thurston.)

Now since $\dim_{\ZZ_2}H_1(M_i);\ZZ_2)\ge6$ for each $i$, we have in
particular that $\Hg(M_i)\ge6$ for each $i$. Since each $M_i$ contains
the genus--$2$ connected closed incompressible surface $S$,
it now follows from Theorem \ref{manganese} that $\vol M_i > 6.45$ for
each $i$. Hence $\vol M > 6.45$.
\EndProof

\Theorem\label{volume3.08}
\threeoeight
\EndTheorem

\Proof
We recall that a group $\Gamma$ is said to be {\it $k$-free} for a
given positive integer $k$ if every subgroup of $\Gamma$ having rank
at most $k$ is free. According to Corollary
\threefreevolume\ of \cite{last}, which was deduced from results in
\cite{ACCS}, if
$M$ is a closed, orientable hyperbolic 3--manifold
such that $\pi_1(M)$ is $3$-free then $\vol M>3.08$. 

Now suppose that $M$ satisfies the hypotheses of Theorem
\ref{volume3.08}, but that $\pi_1(M)$ is not $3$-free. Fix a base
point $P\in M$ and a subgroup $X$ of $\pi_1(M,P)$ which has rank at
most $3$ and is not free. Let $\bar X$ denote the image of $X$ under
the natural homomorphism $\eta:\pi_1(M,P)\to H_1(M;\ZZ_2)$. Then the
subspace $\bar X$ of $H_1(M;\ZZ_2)$ has dimension at most $3$. Since
$\mathrm{dim}_{\mathbb{Z}_2} H_1(M; \mathbb{Z}_2) \geq 4$, there is a
codimension-$1$ subspace $V$ of $H_1(M;\ZZ_2)$ containing $\bar X$.
Then $Y\doteq\eta^{-1}(V)$ is an index-$2$ subgroup of $\pi_1(M,P)$
containing $X$. Hence $Y$ defines a $2$-sheeted based covering space
$p:(\tM,\tP)\to(M,P)$ such that $p_\sharp:\pi_1(\tM,\tP)\to\pi_1((M,P)$
maps some subgroup $\tX$ of $\pi_1(\tM,\tP)$ isomorphically onto
$X$. In particular $\tX$ has rank at most $3$ and is not free, and so
$\pi_1(\tM)$ is not $3$-free.

We now invoke Proposition \homologystuff\ of \cite{CS_onecusp}, which
asserts that if $M$ is a closed, aspherical $3$-manifold, if
$r=\dim_{\ZZ_2}H_1( M;\ZZ_2)$, and if $t$ denotes the dimension of the
image of the cup product map $H^1 (M; \mathbb{Z}_2)\otimes H^1 (M;
\mathbb{Z}_2)\to H^2 (M; \mathbb{Z}_2)$, then for any integer $m\ge0$
and any regular covering $\tM$ of $M$ with covering group $(\ZZ_2)^m$,
we have $\dim_{\ZZ_2}H_1( \tM;\ZZ_2)\ge(m+1)r-m(m+1)/2-t$. Taking $M$
and $\tM$ as above, the hypotheses of of \cite[Proposition
\homologystuff]{CS_onecusp} hold with $m=1$, and by the hypothesis of
the present theorem we have $r\ge4$ and $t\le1$. Hence
$\dim_{\ZZ_2}H_1(\tM;\ZZ_2)\ge6$.

We next invoke Proposition \ifnotwhynot\ of \cite{CS_vol}, which
implies that if $k\ge3$ is an integer and if $N$ is a closed simple
$3$-manifold such that $\dim_{\ZZ_2}H_1( N;\ZZ_2)\ge\max(3k-4,6)$,
then either $\pi_1(N)$ is $k$-free, or $N$ contains a closed connected
incompressible surface of genus at most $k-1$. We may apply this with
$N=\tM$ and $k=3$, since we have seen that
$\dim_{\ZZ_2}H_1(\tM;\ZZ_2)\ge6$. Since we have also seen that
$\pi_1(\tM)$ is not $3$-free, $\tM$ must contain a closed
connected incompressible surface $S$ of genus at most $2$,
and in view of simplicity, $S$ must have genus exactly $2$. Again
using that $\dim_{\ZZ_2}H_1(\tM;\ZZ_2)\ge6$---so that in particular
$\Hg(\tM)\ge6$---we deduce from Theorem \ref{manganese}, with $\tM$
playing the role of $M$, that $\vol \tM>6.45$.  Hence
$$\vol M=\frac12\vol\tM>3.225>3.08.$$
\EndProof

\bibliographystyle{plain}

\end{document}